\documentclass[a4paper, oneside]{amsart}

\RequirePackage{amsmath}
\RequirePackage{bm}
\RequirePackage{amssymb}
\RequirePackage{upref}
\RequirePackage{amsthm}
\RequirePackage{enumerate}
\RequirePackage{pb-diagram}
\RequirePackage{amsfonts}
\RequirePackage[mathscr]{eucal}
\RequirePackage{verbatim}
\RequirePackage{xr}
\RequirePackage{graphicx}
\RequirePackage[pdftex]{floatflt}
\usepackage{calc}
\usepackage{xspace}

\newcommand{\cf}{cf.\@\xspace}
\newcommand{\resp}{resp.\@\xspace}

\newcommand{\aev}{a.e.\@\xspace}

\newcommand{\al}{\alpha}
\newcommand{\bet}{\beta}
\newcommand{\ga}{\gamma}
\newcommand{\de}{\delta }
\newcommand{\e}{\epsilon}

\newcommand{\f}{\varphi}
\newcommand{\h}{\eta}

\newcommand{\ka}{\kappa}
\newcommand{\lam}{\lambda}
\newcommand{\m}{\mu}
\newcommand{\n}{\nu}

\newcommand{\s}{\sigma}
\newcommand{\x}{\xi}

\newcommand{\D}{\varDelta}
\newcommand{\F}{\varPhi}
\newcommand{\Lam}{\varLambda}
\newcommand{\Om}{\varOmega}

\newcommand{\di}[1]{#1\nobreakdash-\hspace{0pt}dimensional}%\di n

\newcommand{\fv}[2]{#1\hspace{0pt}_{|_{#2}}}

\newcommand{\so}{{\mc S_0}}

\newcommand{\const}{\tup{const}}
\newcommand{\ndash}{\nobreakdash--}

\newcommand{\msp[1]}[1]{\mspace{#1mu}}

\newcommand{\R}[1][n+1]{{\protect\mathbb R}^{#1}}

\newcommand{\N}{{\protect\mathbb N}}

\newcommand{\eR}{\stackrel{\lower1ex \hbox{\rule{6.5pt}{0.5pt}}}{\msp[3]\R[]}}
\newcommand{\eN}{\stackrel{\lower1ex \hbox{\rule{6.5pt}{0.5pt}}}{\msp[1]\N}}
\newcommand{\eO}{\stackrel{\lower1ex
\hbox{\rule{6pt}{0.5pt}}}{\msc O}}

\DeclareMathOperator{\graph}{graph}

\newcommand\ra{\rightarrow}

\newcommand\pde[2]{\frac {\partial#1}{\partial#2}}
\newcommand{\un}{\infty}
\newcommand{\A}{\forall}

\newcommand{\set}[2]{\{\,#1\colon #2\,\}}
\newcommand{\uu}{\cup}

\newcommand{\uuu}{\bigcup}

\newcommand{\uud}{ \stackrel{\lower 1ex \hbox {.}}{\uu}}
\newcommand{\uuud}[1]{ \stackrel{\lower 1ex \hbox {.}}{\uuu_{#1}}}
\newcommand\su{\subset}

\newcommand{\sminus}[1][28]{\raise 0.#1ex\hbox{$\scriptstyle\setminus$}}

\newcommand\ti{\times }

\newcommand{\abs}[1]{\lvert#1\rvert}

\newcommand{\norm}[1]{\lVert#1\rVert}

\newcommand{\nnorm}[1]{| \mspace{-2mu} |\mspace{-2mu}|#1| \mspace{-2mu}
|\mspace{-2mu}|}
\newcommand{\spd}[2]{\protect\langle #1,#2\protect\rangle}

\newcommand\ch[3]{\varGamma_{#1#2}^#3}
\newcommand\cha[3]{{\bar\varGamma}_{#1#2}^#3}

\newcommand{\riem}[4]{R_{#1#2#3#4}}
\newcommand{\riema}[4]{{\bar R}_{#1#2#3#4}}

\newcommand{\tit}{\textit}

\newcommand{\tup}{\textup}% text upright

\newcommand{\mc}{\protect\mathcal}
\newcommand{\msc}{\protect\mathscr}

\providecommand{\bysame}{\makebox[3em]{\hrulefill}\thinspace}

\newcommand{\ci}{\cite}

\newcommand{\bt}{\begin{thm}}
\newcommand{\bl}{\begin{lem}}
\newcommand{\bc}{\begin{cor}}
\newcommand{\bd}{\begin{definition}}
\newcommand{\bpp}{\begin{prop}}
\newcommand{\br}{\begin{rem}}
\newcommand{\bn}{\begin{note}}
\newcommand{\be}{\begin{ex}}
\newcommand{\bes}{\begin{exs}}
\newcommand{\bb}{\begin{example}}
\newcommand{\bbs}{\begin{examples}}
\newcommand{\ba}{\begin{axiom}}

\newcommand{\et}{\end{thm}}
\newcommand{\el}{\end{lem}}
\newcommand{\ec}{\end{cor}}
\newcommand{\ed}{\end{definition}}
\newcommand{\epp}{\end{prop}}
\newcommand{\er}{\end{rem}}
\newcommand{\en}{\end{note}}
\newcommand{\ee}{\end{ex}}
\newcommand{\ees}{\end{exs}}
\newcommand{\eb}{\end{example}}
\newcommand{\ebs}{\end{examples}}
\newcommand{\ea}{\end{axiom}}

\newcommand{\bp}{\begin{proof}}
\newcommand{\ep}{\end{proof}}
\newcommand{\eps}{\renewcommand{\qed}{}\end{proof}}

\newcommand{\bal}{\begin{align}}

\newcommand{\bi}[1][1.]{\begin{enumerate}[\upshape #1]}
\newcommand{\bia}[1][(1)]{\begin{enumerate}[\upshape #1]}
\newcommand{\bin}[1][1]{\begin{enumerate}[\upshape\bfseries #1]}
\newcommand{\bir}[1][(i)]{\begin{enumerate}[\upshape #1]}
\newcommand{\bic}[1][(i)]{\begin{enumerate}[\upshape\hspace{2\cma}#1]}
\newcommand{\bis}[2][1.]{\begin{enumerate}[\upshape\hspace{#2\parindent}#1]}
\newcommand{\ei}{\end{enumerate}}

\newcommand\ndots{\raise 0.47ex \hbox {,}\hskip0.06em\cdots %
     \raise 0.47ex \hbox {,}\hskip0.06em} 

\newcommand{\q}{\quad}
\newcommand{\qq}{\qquad}

\newcommand\nd{\noindent}

\newskip\Csmallskipamount                                                
\Csmallskipamount=\smallskipamount
\newskip\Cmedskipamount
\Cmedskipamount=\medskipamount
\newskip\Cbigskipamount
\Cbigskipamount=\bigskipamount

\newcommand\cvs{\vspace\Csmallskipamount}   
\newcommand\cvm{\vspace\Cmedskipamount}

\newskip\csa
\csa=\smallskipamount

\newskip\cma
\cma=\medskipamount

\newskip\cba
\cba=\bigskipamount

\newdimen\spt
\newcommand\citem{\cvs\advance\itemno by
1{(\romannumeral\the\itemno})\hskip3pt}
\newcommand{\bitem}{\cvm\nd\advance\itemno by
1{\bf\the\itemno}\hspace{\cma}}

\newcount\itemno
\newcommand{\las}[1]{\label{S:#1}}

\newcommand{\lae}[1]{\label{E:#1}}

\newcommand{\lal}[1]{\label{L:#1}}
\newcommand{\lad}[1]{\label{D:#1}}

\newcommand{\lar}[1]{\label{R:#1}}

\newcommand{\rs}[1]{Section~\ref{S:#1}}

\newcommand{\rl}[1]{Lemma~\ref{L:#1}}
\newcommand{\rd}[1]{Definition~\ref{D:#1}}

\newcommand{\rr}[1]{Remark~\ref{R:#1}}

\newcommand{\re}[1]{\eqref{E:#1}}

\newskip\thmskip
\thmskip=\parindent

\newskip\hsk
\newenvironment{hinw}{\labelsep=0pt\begin{list}{}{\labelsep=0pt\itemindent=0pt\labelwidth=0pt\leftmargin=\parindent\rightmargin=0pt\partopsep=\cba}%
\item\it\nopagebreak\nopagebreak}%
{\end{list}}

\newcommand\bh{\begin{hinw}}
\newcommand{\eh}{\end{hinw}}

\newtheoremstyle{normal}% name
  {\cba}%      Space above, empty = `usual value'
  {\cba}%      Space below
  {}% Body font
  {\thmskip}%Indent amount (empty = no indent, \parindent = para indent)
  {\bfseries}% Thm head font
  {.}%        Punctuation after thm head
  {\hsk}%     Space after thm head: " " = normal interword space;
  {}% Thm head spec

\newtheoremstyle{abschnitt}% name
  {\cba}%      Space above, empty = `usual value'
  {\cba}%      Space below
  {}% Body font
  {\thmskip}% Indent amount (empty = no indent, \parindent = para indent)
  {\bfseries}% Thm head font
  {.}%        Punctuation after thm head
  {\hsk}%     Space after thm head: " " = normal interword space;
  {}% Thm head spec

\newtheoremstyle{italic}% name
  {\cba}%      Space above, empty = `usual value'
  {\cba}%      Space below
  {\itshape}% Body font
  {\thmskip}%  Indent amount (empty = no indent, \parindent = para indent)
  {\bfseries}% Thm head font
  {.}%        Punctuation after thm head
  {\hsk}%     Space after thm head: " " = normal interword space;
  {}% Thm head spec

\newtheoremstyle{aufgaben}% name
  {\cba}%      Space above, empty = `usual value'
  {\cba}%      Space below
  {}% Body font
  {}%         Indent amount (empty = no indent, \parindent = para indent)
  {\normalsize\bfseries}% Thm head font
  {.}%        Punctuation after thm head
  {\hsk}%     Space after thm head: " " = normal interword space;
  {}% Thm head spec

\newtheoremstyle{break}% name
  {\cba}%      Space above, empty = `usual value'
  {\cba}%      Space below
  {\itshape}% Body font
  {}%         Indent amount (empty = no indent, \parindent = para indent)
  {\bfseries}% Thm head font
  {.}%        Punctuation after thm head
  {\newline}% Space after thm head: \newline = linebreak
  {}%         Thm head spec

\swapnumbers
\theoremstyle{italic}
\newtheorem{thm}[subsection]{Theorem}
\newtheorem{lem}[subsection]{Lemma}
\newtheorem{prop}[subsection]{Proposition}
\newtheorem{cor}[subsection]{Corollary}

\theoremstyle{normal}
\newtheorem{rem}[subsection]{Remark}
\newtheorem{definition}[subsection]{Definition}
\newtheorem{example}[subsection]{Example}
\newtheorem{examples}[subsection]{Examples}
\newtheorem{ex}[subsection]{Exercise}
\newtheorem{note}[subsection]{}
\newtheorem{axiom}[subsection]{Axiom}

\theoremstyle{aufgaben}
\newtheorem{exs}[subsection]{Exercises}

\swapnumbers

\numberwithin{equation}{section}
\numberwithin{figure}{section}

\newenvironment{textequation}[1][0.8]
{\begin{equation}
\begin{aligned}
\begin{minipage}{#1\linewidth}}
{\end{minipage}
\end{aligned}
\end{equation}
\ignorespacesafterend}

\newcommand{\btext}{\begin{textequation}}
\newcommand{\etext}{\end{textequation}}

\usepackage[german,english]{babel}
\usepackage{graphicx}
\RequirePackage{amsmath}
\RequirePackage{bm}
\RequirePackage{amssymb}
\RequirePackage{upref}
\RequirePackage{amsthm}
\RequirePackage{enumerate}
\RequirePackage{pb-diagram}
\RequirePackage{amsfonts}
\RequirePackage[mathscr]{eucal}
\RequirePackage{verbatim}
\RequirePackage{xr}
\RequirePackage{graphicx}
\RequirePackage[pdftex]{floatflt}
\usepackage{calc}
\usepackage{xspace}
\RequirePackage{upref}
\RequirePackage{amsthm}
\RequirePackage{enumerate}%\begin{enumerate}[(i)]
\usepackage[mathscr]{eucal}
\usepackage{xr-hyper}

\usepackage{calc}

\newlength{\oddsidemarginlength}
\newlength{\topmarginlength}

\newcounter{numberoflines}
\newcounter{tempcc}
\oddsidemargin=\oddsidemarginlength
\evensidemargin=\oddsidemargin
\topmargin=\topmarginlength
\begin{document}

\flushbottom

%\larger[1]
%\frontmatter

\title{The inverse mean curvature flow in cosmological spacetimes}

% author one information
\author{Claus Gerhardt}
\address{Ruprecht-Karls-Universit\"at, Institut f\"ur Angewandte Mathematik,
Im Neuenheimer Feld 294, 69120 Heidelberg, Germany}
%\curraddr{}
\email{gerhardt@math.uni-heidelberg.de}
\urladdr{http://www.math.uni-heidelberg.de/studinfo/gerhardt/}
\thanks{This work has been supported by the Deutsche Forschungsgemeinschaft.}

% author two information
%\author{}
%\address{}
%\curraddr{}
%\email{}
%\thanks{}
%
\date{\today}
%
% at present the "communicated by" line appears only in ERA and PROC
%\commby{}

%\dedicatory{}
\begin{abstract}
We prove that the leaves of an inverse mean curvature flow provide a foliation of a future end of a cosmological spacetime $N$ under the necessary and sufficent assumptions that $N$ satisfies a future mean curvature barrier condition and a strong volume decay condition. Moreover, the flow parameter $t$ can be used to define a new physically important time function.
\end{abstract}
\maketitle

%\include{schmutztitelI}
%\include{titleI}
%\include{impressumI}
%\include{preface}
%\maketitle

\tableofcontents

\section{Introduction}

The inverse mean curvature flow has already been considered in Euclidean space
\cite{cg90} or in asymptotically flat Riemannian spaces \cite{hi:penrose}. In the latter
case Huisken and Ilmanen used it to prove the  Penrose inequality. One major
difficulty in their proof was that jumps might occur during the flow, i.e., the mean
curvature of the flow hypersurfaces might vanish even though the initial
hypersurface has positive mean curvature.

\cvm
The Lorentzian geometry is much more favourable for curvature flows, \cf \cite{eh1,
cg:indiana, cg:mz, cg:scalar}, so that no jumps should occur in case of the inverse
mean curvature flow. We shall show that this is indeed the case, if the ambient space is
a \tit{globally hyperbolic} $(n+1)$-dimensional Lorentzian manifold $N$ with a
compact Cauchy hypersurface satisfying the timelike convergence condition
\begin{equation}
\bar R_{\al\bet}\nu^\al\nu^\bet \ge 0\qq\A\;\spd\nu\nu=-1.
\end{equation}

Such spaces are called \tit{cosmological spacetimes}, a terminology due to Bartnik.

\cvm
Let $M_0\su N$ be a spacelike hypersurface the mean curvature of which is either
strictly positive or negative, then we consider the inverse mean curvature flow (IMCF)
\begin{equation}\lae{0.2}
\dot x=-H^{-1}\nu
\end{equation}
with initial hypersurface $M_0$. Here, $\nu$ is the past directed normal of the flow
hypersurfaces $M(t)$ and $H=\fv H{M(t)}$ the corresponding mean curvature, i.e.,
the trace of the second fundamental form.

If $\fv H{M_0}$ is positive \resp negative, then the flow moves to the future \resp
past of $M_0$. Furthermore, $\fv H{M(t)}$ will uniformly tend to $\un$ \resp
$-\un$, if the flow exists for all time.

In former papers we referred to this latter phenomenen by saying that there were
\tit{crushing singularities} in the future  \resp past, erroneously assuming that only
big crunch or big bang type singularities could produce spacelike hypersurfaces the
mean curvatures of which become unbounded if the hypersurfaces approached the
singularities.

But a behaviour like that could also be caused by a \tit{null hypersurface} $\mc H$,
e.g., by the event horizon of a black hole, if the spacetime can be viewed as having a
past or future boundary component $\mc H$ that can be identified with a compact
null hypersurface representing a \tit{non-crushing} singularity, i.e., the Riemannian
curvature tensors remains uniformly bounded near $\mc H$
\begin{equation}
\riema\al\bet\ga\de \bar R^{\al\bet\ga\de}\le \const.
\end{equation}
An example of such a spacetime is given in \rs{1}.

\cvm
We therefore define

\bd
Let $N$ be a globally hyperbolic spacetime with compact Cauchy hypersurface $\so$
so that $N$ can be written as a topological product $N=\R[]\times \so$ and its
metric expressed as
\begin{equation}\lae{0.4}
d\bar s^2=e^{2\psi}(-(dx^0)^2+\s_{ij}(x^0,x)dx^idx^j).
\end{equation}
Here, $x^0$ is a globally defined future directed time function and $(x^i)$ are local
coordinates for $\so$.
$N$ is said to have a \tit{future mean curvature barrier} \resp \tit{past mean
curvature barrier}, if there are sequences $M_k^+$ \resp $M_k^-$ of closed
spacelike hypersurfaces such that
\begin{equation}
\lim_{k\ra\un} \fv H{M_k^+}=\un \q\tup{\resp}\q \lim_{k\ra\un} \fv H{M_k^-}=-\un
\end{equation}
and
\begin{equation}
\limsup \inf_{M_k^+}x^0> x^0(p)\qq\A\,p\in N
\end{equation}
\resp
\begin{equation}
\liminf \sup_{M_k^-}x^0< x^0(p)\qq\A\,p\in N.
\end{equation}
\ed

\br\lar{0.2}
Let $N$ be a cosmological spacetime with  future and a past mean curvature barriers,
then it can be foliated by closed hypersurfaces of constant mean curvature, \cf
\cite{cg1}. Moreover, the mean curvature function $\tau$ is continuous in $N$ and
smooth in $\{\tau\ne 0\}$ with non-vanishing gradient, hence it can be used as a
time function, \cf \cite{cg:foliation}. These results are also valid in future \resp past
ends. 
\er

We shall assume in the following that $N$ has a future mean curvature barrier. By
reversing the time direction this configuration also comprises the case that $N$ has a
past mean curvature barrier.

Under this assumption we shall prove that, for a given compact spacelike
hypersurface $M_0$ with $\fv H{M_0}>0$, the future of $M_0$ can be foliated by
the leaves of an IMCF starting at $M_0$ provided a so-called future \tit{strong volume
decay condition} is satisfied, \cf \rd{1.2}. A strong volume decay condition is both
necessary and sufficient in order that the IMCF exists for all time.

\cvm
The main result of this paper can be summarized in the following theorem

\bt
Let $N$ be a cosmological spacetime with compact Cauchy hypersurface $\so$ and
with a  future mean curvature barrier. Let $M_0$ be a closed spacelike hypersurface
with positive mean curvature and assume furthermore that $N$ satisfies a future
volume decay condition. Then the IMCF \re{0.2} with initial hypersurface $M_0$ exists
for all time and provides a foliation of the future $D^+(M_0)$ of $M_0$.

The evolution parameter $t$ can be chosen as a new time function. The flow
hypersurfaces $M(t)$ are the slices $\{t=\const\}$ and their volume satisfies
\begin{equation}
\abs{M(t)}=\abs{M_0} e^{-t}.
\end{equation}

Defining an almost proper time function $\tau$ by choosing
\begin{equation}
\tau=1-e^{-\frac1n t}
\end{equation}
we obtain  $0\le \tau <1$,
\begin{equation}
\abs{M(\tau)}=\abs{M_0} (1-\tau)^n,
\end{equation}
and  the future singularity corresponds to $\tau=1$.

Moreover, the length $L(\ga)$ of any future directed curve $\ga$ starting from
$M(\tau)$ is bounded from above by
\begin{equation}
L(\ga)\le c (1-\tau),
\end{equation}
where $c=c(n, M_0)$. Thus, the expression $1-\tau$ can be looked at as the radius
of the slices $\{\tau=\const\}$ as well as a measure of the remaining life span of the
universe.
\et

Without any further structural assumptions it seems impossible to derive any
convergence results for an appropriately rescaled IMCF. In \cite{cg:arw} we look at
the IMCF in asymptotically Robertson Walker spaces and prove that a properly rescaled
flow converges indeed.

\section{Notations and definitions}\las{1}

The main objective of this section is to state the equations of Gau{\ss}, Codazzi,
and Weingarten for spacelike hypersurfaces $M$ in a \di {(n+1)} Lorentzian
manifold
$N$.  Geometric quantities in $N$ will be denoted by
$(\bar g_{ \al \bet}),(\riema  \al \bet \ga \de)$, etc., and those in $M$ by $(g_{ij}), 
(\riem ijkl)$, etc.. Greek indices range from $0$ to $n$ and Latin from $1$ to $n$;
the summation convention is always used. Generic coordinate systems in $N$ resp.
$M$ will be denoted by $(x^ \al)$ \resp $(\x^i)$. Covariant differentiation will
simply be indicated by indices, only in case of possible ambiguity they will be
preceded by a semicolon, i.e., for a function $u$ in $N$, $(u_ \al)$ will be the
gradient and
$(u_{ \al \bet})$ the Hessian, but e.g., the covariant derivative of the curvature
tensor will be abbreviated by $\riema  \al \bet \ga{ \de;\e}$. We also point out that
\begin{equation}
\riema  \al \bet \ga{ \de;i}=\riema  \al \bet \ga{ \de;\e}x_i^\e
\end{equation}
with obvious generalizations to other quantities.

Let $M$ be a \tit{spacelike} hypersurface, i.e., the induced metric is Riemannian,
with a differentiable normal $\n$ which is timelike.

In local coordinates, $(x^ \al)$ and $(\x^i)$, the geometric quantities of the
spacelike hypersurface $M$ are connected through the following equations
\begin{equation}\lae{1.2}
x_{ij}^ \al= h_{ij}\n^ \al
\end{equation}
the so-called \tit{Gau{\ss} formula}. Here, and also in the sequel, a covariant
derivative is always a \tit{full} tensor, i.e.

\begin{equation}
x_{ij}^ \al=x_{,ij}^ \al-\ch ijk x_k^ \al+ \cha  \bet \ga \al x_i^ \bet x_j^ \ga.
\end{equation}
The comma indicates ordinary partial derivatives.

In this implicit definition the \tit{second fundamental form} $(h_{ij})$ is taken
with respect to $\n$.

The second equation is the \tit{Weingarten equation}
\begin{equation}
\n_i^ \al=h_i^k x_k^ \al,
\end{equation}
where we remember that $\n_i^ \al$ is a full tensor.

Finally, we have the \tit{Codazzi equation}
\begin{equation}
h_{ij;k}-h_{ik;j}=\riema \al \bet \ga \de\n^ \al x_i^ \bet x_j^ \ga x_k^ \de
\end{equation}
and the \tit{Gau{\ss} equation}
\begin{equation}
\riem ijkl=- \{h_{ik}h_{jl}-h_{il}h_{jk}\} + \riema  \al \bet\ga \de x_i^ \al x_j^ \bet
x_k^ \ga x_l^ \de.
\end{equation}

Now, let us assume that $N$ is a globally hyperbolic Lorentzian manifold with a
\tit{compact} Cauchy surface. 
$N$ is then a topological product $I\times \mc S_0$, where $I$ is an open interval,
$\mc S_0$ is a compact Riemannian manifold, and there exists a Gaussian coordinate
system
$(x^ \al)$, such that the metric in $N$ has the form 
\begin{equation}\lae{1.7}
d\bar s_N^2=e^{2\psi}\{-{dx^0}^2+\s_{ij}(x^0,x)dx^idx^j\},
\end{equation}
where $\s_{ij}$ is a Riemannian metric, $\psi$ a function on $N$, and $x$ an
abbreviation for the spacelike components $(x^i)$. 
We also assume that
the coordinate system is \tit{future oriented}, i.e., the time coordinate $x^0$
increases on future directed curves. Hence, the \tit{contravariant} timelike
vector $(\x^ \al)=(1,0,\dotsc,0)$ is future directed as is its \tit{covariant} version
$(\x_ \al)=e^{2\psi}(-1,0,\dotsc,0)$.

Let $M=\graph \fv u\so$ be a spacelike hypersurface
\begin{equation}
M=\set{(x^0,x)}{x^0=u(x),\,x\in\mc S_0},
\end{equation}
then the induced metric has the form
\begin{equation}
g_{ij}=e^{2\psi}\{-u_iu_j+\s_{ij}\}
\end{equation}
where $\s_{ij}$ is evaluated at $(u,x)$, and its inverse $(g^{ij})=(g_{ij})^{-1}$ can
be expressed as
\begin{equation}\lae{1.10}
g^{ij}=e^{-2\psi}\{\s^{ij}+\frac{u^i}{v}\frac{u^j}{v}\},
\end{equation}
where $(\s^{ij})=(\s_{ij})^{-1}$ and
\begin{equation}\lae{1.11}
\begin{aligned}
u^i&=\s^{ij}u_j\\
v^2&=1-\s^{ij}u_iu_j\equiv 1-\abs{Du}^2.
\end{aligned}
\end{equation}
Hence, $\graph u$ is spacelike if and only if $\abs{Du}<1$.

The covariant form of a normal vector of a graph looks like
\begin{equation}
(\n_ \al)=\pm v^{-1}e^{\psi}(1, -u_i).
\end{equation}
and the contravariant version is
\begin{equation}
(\n^ \al)=\mp v^{-1}e^{-\psi}(1, u^i).
\end{equation}
Thus, we have
\br Let $M$ be spacelike graph in a future oriented coordinate system. Then the
contravariant future directed normal vector has the form
\begin{equation}
(\n^ \al)=v^{-1}e^{-\psi}(1, u^i)
\end{equation}
and the past directed
\begin{equation}\lae{1.15}
(\n^ \al)=-v^{-1}e^{-\psi}(1, u^i).
\end{equation}
\er

In the Gau{\ss} formula \re{1.2} we are free to choose the future or past directed
normal, but we stipulate that we always use the past directed normal for reasons
that we have explained in \ci[Section 2]{cg:indiana}.

Look at the component $ \al=0$ in \re{1.2} and obtain in view of \re{1.15}

\begin{equation}\lae{1.16}
e^{-\psi}v^{-1}h_{ij}=-u_{ij}- \cha 000\mspace{1mu}u_iu_j- \cha 0j0
\mspace{1mu}u_i- \cha 0i0\mspace{1mu}u_j- \cha ij0.
\end{equation}
Here, the covariant derivatives are taken with respect to the induced metric of
$M$, and
\begin{equation}
-\cha ij0=e^{-\psi}\bar h_{ij},
\end{equation}
where $(\bar h_{ij})$ is the second fundamental form of the hypersurfaces
$\{x^0=\const\}$.

An easy calculation shows
\begin{equation}
\bar h_{ij}e^{-\psi}=-\tfrac{1}{2}\dot\s_{ij} -\dot\psi\s_{ij},
\end{equation}
where the dot indicates differentiation with respect to $x^0$.

\cvm
Next we shall define the strong volume decay condition.

\bd\lad{1.2}
Suppose there exists a time function $x^0$ such that the future end of $N$ is
determined by $\{\tau_0\le x^0<b\}$ and the coordinate slices
$M_\tau=\{x^0=\tau\}$ have positive mean curvature with respect to the past
directed normal for $\tau_0\le\tau<b$. In addition the volume $\abs{M_\tau}$
should satisfy
\begin{equation}\lae{1.19}
\lim_{\tau\ra b}\abs{M_\tau}=0.
\end{equation}

A decay like that is normally associated with a future singularity and we simply call it
\tit{volume decay}. If $(g_{ij})$ is the induced metric of $M_\tau$ and
$g=\det(g_{ij})$, then we have
\begin{equation}\lae{1.20}
\log g(\tau_0,x)-\log g(\tau,x)=\int_{\tau_0}^\tau e^\psi \bar H(s,x)\q\A\,x\in \so,
\end{equation}
where $\bar H(\tau,x)$ is the mean curvature of $M_\tau$ in $(\tau,x)$. For a proof
we refer to \cite{cg:volume}. 

\cvm
In view of \re{1.19} the left-hand side of this equation tends to infinity if $\tau$
approaches $b$ for \aev $x\in \so$, i.e.,
\begin{equation}
\lim_{\tau\ra b}\int_{\tau_0}^\tau e^\psi \bar H(s,x)=\un\q \tup{for
\aev}\; x\in\so.
\end{equation}

Assume now, there exists a continuous, positive function $\f=\f(\tau)$ such that
\begin{equation}\lae{1.22}
e^\psi \bar H(\tau,x)\ge \f(\tau)\qq\A\, (\tau, x)\in (\tau_0,b)\times \so,
\end{equation}
where
\begin{equation}\lae{1.23}
\int_{\tau_0}^b \f(\tau)=\un,
\end{equation}
then we say that the future of $N$ satisfies a \tit{strong volume decay condition}.
\ed

\br
(i) By approximation we may---and shall---assume that the function $\f$ above is
smooth.

(ii) A similar definition holds for the past of $N$ by simply reversing the time
direction. Notice that in this case the mean curvature of the coordinate slices has to
be negative.
\er

\bl\lal{1.4}
Suppose that the future of $N$ satisfies a strong volume decay condition, then there
exist a time function $\tilde x^0=\tilde x^0(x^0)$, where $x^0$ is the time function
in the strong volume decay condition, such that the mean curvature $\bar H$ of the
slices $\tilde x^0=\const$ satisfies the estimate
\begin{equation}\lae{1.24}
e^{\tilde\psi}\bar H\ge 1.
\end{equation}
The factor $e^{\tilde\psi}$ is now the conformal factor in the representation
\begin{equation}\lae{1.25}
d\bar s^2=e^{2\tilde\psi}(-(d\tilde x^0)^2+\s_{ij}dx^idx^j).
\end{equation}

The range of $\tilde x^0$ is equal to the interval $[0,\un)$, i.e., the singularity
corresponds to $\tilde x^0=\un$.
\el

\bp
Define $\tilde x^0$ by
\begin{equation}
\tilde x^0=\int_{\tau_0}^{x^0}\f(\tau),
\end{equation}
where $\f$ is the function in \re{1.22} now assumed to be smooth.

The conformal factor in \re{1.25} is then equal to
\begin{equation}
e^{2\tilde\psi}=e^{2\psi} \pde{x^0}{\tilde x^0} \pde{x^0}{\tilde x^0}=e^{2\psi}
\f^{-2},
\end{equation}
and hence
\begin{equation}
e^{\tilde\psi}\bar H=e^\psi \bar H \f^{-1}\ge 1,
\end{equation}
in view of \re{1.22}.
\ep

As we mentioned in the introduction there are spacetimes which satisfy a mean
curvature barrier condition but the resulting singularity is not crushing.

\cvm
To construct an example let us start with a  $\tup{S-AdS}_{(n+2)}$ spacetime with
metric
\begin{equation}
d\hat s^2=-f dt^2+ f^{-1}dr^2+r^2 \s_{ij}dx^idx^j,
\end{equation}
where
\begin{equation}
f=\ka-\frac2{n(n+1)}\Lam r^2-m r^{-(n-1)}
\end{equation}
with constants $\Lam$ and $m>0$; $(\s_{ij})$ is the metric of a compact
$n$-dimensional spaceform of curvature $\ka=0,1,-1$.

This spacetime satisfies the Einstein equations
\begin{equation}
G_{\al\bet}+\Lam \bar g_{\al\bet}=0.
\end{equation}

\cvm
Let us suppose for simplicity that $\ka=1$ and $\Lam<0$, though this is not
important in  our considerations. In $\{r=0\}$ is a black hole singularity and the event
horizon $\mc H=f^{-1}(0)$ is characterized by $r=r_0$. 

The region $\{f<0\}$ is the black hole region. In this region $r$ is the time function
and $t$ is a spatial variable. Let us pick the black hole region. 

Normally the variable
$t$ describes the real axis, but, since it is a spatial  variable, we are free to
compactify it, and we shall suppose that $t$ is a variable for $S^1$. By this
compactification we have defined a globally hyperbolic spacetime $N$ with compact
Cauchy hypersurface $\so=S^1\times S^n$ which satisfies the timelike convergence
condition since
\begin{equation}
\bar R_{\al\bet}=\tfrac2n \Lam \bar g_{\al\bet}
\end{equation}
and $\Lam$ is supposed to be negative.

\cvm
$N$ has a crushing singularity in $r=0$, and, as we shall show in a moment, also a
mean curvature barrier singularity in $r=r_0$, which is however not crushing, since
the metric quantities were not changed by the compactifcation but only the topology.

Define
\begin{equation}
\tilde f=-f\q\tup{and}\q \psi=-\tfrac12 \log\tilde f,
\end{equation}
then the metric can be expressed as
\begin{equation}
\begin{aligned}
d\bar s^2&=e^{2\psi}(-dr^2+\tilde f^2dt^2+\tilde f r^2\s_{ij} dx^idx^j)\\
&\equiv e^{2\psi}(-dr^2+\tilde\s_{ab}dx^adx^b).
\end{aligned}
\end{equation}

The second fundamental form of the hypersurfaces $\{r=\const\}$ with respect to
the past directed normal is given by
\begin{equation}
e^{-\psi}\bar h_{ab}=\tfrac12 \dot{\tilde\s}_{ab}-\tfrac12 \tilde f^{-1}\dot{\tilde
f} \tilde\s_{ab},
\end{equation}
where the dot indicates differentiation with respect to $r$, and where we note that
the time function $r$ is past directed in contrast to the usual convention. Hence the
mean curvature $\bar H$ is equal to
\begin{equation}
\bar H=\tilde f^{-\frac12}(\tfrac12 \dot{\tilde f}+n\tilde f r^{-1})
\end{equation}
and we deduce that $\bar H$ tends to $-\un$, if the hypersurfaces approach the
horizon $\mc H$, and to $\un$, if the hypersurfaces approach the black hole
singularity $r=0$.

\cvm
Sometimes, we need a Riemannian reference metric, e.g., if we want to estimate
tensors. Since the Lorentzian metric can be expressed as
\begin{equation}
\bar g_{\al\bet}dx^\al dx^\bet=e^{2\psi}\{-{dx^0}^2+\s_{ij}dx^i dx^j\},
\end{equation}
we define a Riemannian reference metric $(\tilde g_{\al\bet})$ by
\begin{equation}
\tilde g_{\al\bet}dx^\al dx^\bet=e^{2\psi}\{{dx^0}^2+\s_{ij}dx^i dx^j\}
\end{equation}
and we abbreviate the corresponding norm of a vectorfield $\h$ by
\begin{equation}
\nnorm \h=(\tilde g_{\al\bet}\h^\al\h^\bet)^{1/2},
\end{equation}
with similar notations for higher order tensors.

\section{The evolution problem}

The evolution problem \re{0.2} is a parabolic problem, hence a solution exists
 on a maximal time interval $[0,T^*)$, $0<T^*\le \un$, cf.
\ci[Section 2]{cg96}, where we apologize for the ambiguity of also calling the
evolution parameter \tit{time}.

Next, we want to show how the metric, the second fundamental form, and the
normal vector of the hypersurfaces $M(t)$ evolve. All time derivatives are
\tit{total} derivatives. We refer to \ci{cg:indiana} for more general results and to
\ci[Section 3]{cg96}, where proofs are given in a Riemannian setting, but these
proofs are also valid in a Lorentzian environment.

\bl
The metric, the normal vector, and the second fundamental form of $M(t)$
satisfy the evolution equations
\begin{equation}\lae{2.1}
\dot g_{ij}=-2 H^{-1}h_{ij},
\end{equation}
\begin{equation}\lae{2.2}
\dot \n=\nabla_M(-H^{-1})=g^{ij}(-H^{-1})_i x_j,
\end{equation}
and
\begin{equation}\lae{2.3}
\dot h_i^j=(-H^{-1})_i^j+H^{-1} h_i^k h_k^j + H^{-1} \riema
\al\bet\ga\de\n^\al x_i^\bet \n^\ga x_k^\de g^{kj}
\end{equation}
\begin{equation}
\dot h_{ij}=(-H^{-1})_{ij}-H^{-1} h_i^k h_{kj}+ H^{-1} \riema
\al\bet\ga\de\n^\al x_i^\bet \n^\ga x_j^\de.
\end{equation}
\el

\bl[Evolution of $H^{-1}$]\lal{2.2}
The term $H^{-1}$ evolves according to the equation
\begin{equation}\lae{2.5}
\begin{aligned}
{(H^{-1})}^\prime- H^{-2}\D H^{-1}=&\msp[0]- H^{-2}(\norm A^2+\bar
R_{\al\bet}\nu^\al\nu^\bet) H^{-1}
\end{aligned}
\end{equation}
where
\begin{equation}
(H^{-1})^{\prime}=\frac{d}{dt}H^{-1}
\end{equation}
and
\begin{equation}
\norm A^2=h_{ij}h^{ij}.
\end{equation}
\el

From \re{2.3} we deduce with the help of the Ricci identities and the Codazzi
equations a parabolic equation for the second fundamental form

\bl\lal{2.3}
The mixed tensor $h_i^j$ satisfies the parabolic equation
\begin{equation}\lae{2.8}
\begin{aligned}
\dot h_i^j-&H^{-2}\D h_i^j\\
&=-H^{-2}\norm A^2h_i^j+2 H^{-1} h_i^kh_k^j
-2 H^{-3}
H_i H^j\\ 
&\q\,+2H^{-2}\riema \al\bet\ga\de x_m^\al x_i ^\bet x_k^\ga
x_r^\de h^{km} g^{rj}\\
&\q\,-H^{-2}g^{kl}\riema \al\bet\ga\de x_m^\al x_k ^\bet x_r^\ga x_l^\de
h_i^m g^{rj}\\
&\q\,- H^{-2}g^{kl}\riema \al\bet\ga\de x_m^\al x_k ^\bet x_i^\ga
x_l^\de h^{mj} \\ 
&\q\,-H^{-2}\bar R_{\al\bet}\n^\al\n^\bet h_i^j+2 H^{-1}\riema
\al\bet\ga\de\n^\al x_i^\bet\n^\ga x_m^\de g^{mj}\\ 
&\q\,+ H^{-2}g^{kl}\bar R_{\al\bet\ga\de;\e}\{\n^\al x_k^\bet x_l^\ga x_i^\de
x_m^\e g^{mj}+\n^\al x_i^\bet x_k^\ga x_m^\de x_l^\e g^{mj}\}.
\end{aligned}
\end{equation}
\el

Since the timelike convergence condition is assumed to be valid we immediately
deduce from \rl{2.2}

\bl\lal{2.4}
There exists a positive constant $c_0=c_0(M_0)$, such that  the
estimate
\begin{equation}\lae{2.9}
H\ge c_0\, e^{\frac1n t}
\end{equation}
is valid during the evolution.
\el

\bp
Let $\f= H^{-1} e^{\frac1n t}$, then $\f$ satisfies the inequality
\begin{equation}
\dot \f -H^{-2}\D\f\le -H^{-2}\abs{A}^2\f+\tfrac1n\f\le 0,
\end{equation}
hence we conclude
\begin{equation}
\f\le \sup_{M_0}\f=\sup_{M_0}H.\qedhere
\end{equation}
\ep

\section{Lower order estimates}

The evolution problem \re{0.2} exists on a maximal time interval $I=[0,T^*)$. We
want to prove that $T^*=\un$, and that the flow hypersurfaces $M(t)$ run into the
future singularity, if $t$ tends to infinity.

The latter property is a characteristicum of the inverse mean curvature flow under
very weak assumptions: if the flow exists for all time, then it cannot stay in a compact
region of $N$, or, more precisely

\bl\lal{3.1}
Let $N$ be a cosmological spacetime with a future mean curvature barrier, and let
$M_0$ be a compact spacelike hypersurface with positive mean curvature. Suppose
that $N=\R[]\times\so$ and that the metric is given as in \re{0.4}. Assume that
the inverse mean curvature flow with initial hypersurface $M_0$ exists for all time,
and let the flow hypersurfaces $M(t)$ be expressed as graphs of a function $u$ over
$\so$
\begin{equation}
M(t)=\set{(x^0,x)}{x^0=u(t,x),\,x\in \so}.
\end{equation}
Then there holds
\begin{equation}\lae{3.2}
\lim_{t\ra\un} \inf_\so u(t,\cdot)=\un.
\end{equation}
\el

\bp
(i) Because of the barrier condition a future end of $N$, $N_+$, can be foliated by
hypersurfaces of positive constant mean curvature and we can choose the mean curvature
$\tau$ of that CMC foliation as new time function $x^0=\tau$ in $N_+$
\begin{equation}
N_+=\set{(\tau,x)}{k\le \tau<\un,\, x\in\so},
\end{equation}
\cf \rr{0.2}, where $k$ is a positive constant and where we used the same symbol
$\so$ for the compact Cauchy hypersurface---indeed, we could use the original Cauchy
hypersurface $\so$, since it need not be a level hypersurface.

\cvm
Let $t_0$ be such that
\begin{equation}\lae{3.4}
c_0 e^{\frac1n t_0}>2k,
\end{equation}
where $c_0$ is the constant in inequality \re{2.9}, then we claim that
\begin{equation}
M(t)\su N_+\qq\A\,t\ge t_0.
\end{equation}

To prove this claim we shall apply the Synge's lemma.
Denote the coordinate slices  $x^0=\tau$ by $M_\tau$, i.e., $M_\tau$ has
constant mean curvature $\bar H=\tau$. 

\cvm
It suffices to show that all $M(t)$ with
$t\ge t_0$ lie in the future of $M_k$. Suppose this were not the case for some
$M(t)$, then the Lorentzian distance between $M(t)$ and $M_k$ would be positive
\begin{equation}
d=d(M(t),M_k)>0
\end{equation}
and hence there would exist a maximal future directed geodesic $\ga$ from $M(t)$
to $M_k$. Synge's lemma would then yield
\begin{equation}
\fv H{M_k}(\ga (d))\ge \fv H{M(t)}(\ga (0))+\int_0^d \bar R_{\al\bet}\dot\ga^\al
\dot\ga^\bet;
\end{equation}
a contradiction in view of \re{3.4} and the timelike convergence condition.

\cvm
(ii) Thus, the flow hypersurfaces $M(t)$ are covered by the new coordinate system
for $t\ge t_0$. The metric of $N$ has again the form as in \re{0.4}. 

Now, the mean curvature  $\bar H$ of the coordinate slices satisfies the evolution
equation
\begin{equation}\lae{3.8}
\dot{\bar H}=-\D e^\psi +(\abs{\bar A}^2+\bar R_{\al\bet}\nu^\al \nu^\bet)
e^\psi,
\end{equation}
where the dot indicates differentiation with respect to $x^0$, the Laplace operator is
the Laplace Beltrami operator of the slice, $\abs{\bar A}^2$ the square of the second
fundamental form and $\nu$ the past directed normal and $e^\psi$ the conformal
factor of the metric.

 This relation is valid for the
slices of any time function
$x^0$ for which the metric has the form as in
\re{0.4}, since the slices are solutions of the evolution equation
\begin{equation}
\dot x=-e^\psi \nu
\end{equation}
from which the relation \re{3.8} can be easily deduced:  apply  the general
formula (3.8) in \cite{cg:indiana} with $\F=-e^\psi$.

\cvm
For the special time function $x^0=\tau$ we therefore obtain
\begin{equation}
1=\dot{\bar H}\ge -\D e^{\psi}+\tfrac1n \tau^2 e^{\psi}.
\end{equation}

Moreover, let $x_0\in\so$ be a point where, for fixed $\tau$,
\begin{equation}
\sup_\so e^{\psi(\tau,\cdot)}=e^{\psi(\tau,x_0)},
\end{equation}
then the maximum principle implies
\begin{equation}
1\ge \tfrac1n \tau^2 e^{\psi(\tau,x_0)}\ge \tfrac1n \tau^2
e^{\psi(\tau,x)}\qq\A\,x\in\so
\end{equation}
and hence
\begin{equation}\lae{3.13}
\bar He^\psi\le n \bar H^{-1}
\end{equation}
for all slices $M_\tau$.

This inequality will be the key ingredient to prove the limit relation \re{3.2}.

\cvm
(iii) Define the function $\f$ on $t\ge t_0$ by
\begin{equation}
\f(t)=\inf_\so u(t,\cdot),
\end{equation}
then $\f$ is Lipschitz continuous and for \aev $t$ there holds
\begin{equation}
\dot\f(t)=\dot u(t,x_t),
\end{equation}
where $x_t$ is such that the infimum is attained in $x_t$. This result is well known;
we shall give a short prove in \rl{3.2} below for the sake of completeness.

Now, from \re{0.2}, looking at the component $\al=0$, we deduce that $u$ satisfies
the evolution equation
\begin{equation}
\dot u=\frac{\tilde v}{H e^\psi},
\end{equation}
where $\tilde v=v^{-1}$ and where the time derivative is the total derivative, i.e.,
\begin{equation}
\dot u=\pde ut +u_i\dot x^i
\end{equation}
and hence
\begin{equation}\lae{3.18}
\pde ut =\frac v{H e^\psi}.
\end{equation}

From \re{1.16} we infer
\begin{equation}
e^{-\psi}\tilde v H=-\D u-\cha 000\norm{Du}^2-2 \cha 0i0 u^i +e^{-\psi}\bar H,
\end{equation}
and conclude further, with the help of the maximum principle, that in $x_t$
\begin{equation}
H\le \bar H,
\end{equation}
and thus
\begin{equation}
\pde ut\ge \frac1{\bar H e^\psi}
\end{equation}
in $x_t$.

Therefore, $\f$ satisfies
\begin{equation}
\dot\f\ge \frac1{\bar H e^\psi}\qq\tup{for \aev}\; t\ge t_0,
\end{equation}
hence
\begin{equation}\lae{3.23}
\dot\f\ge \tfrac1n \bar H=\tfrac1n \f
\end{equation}
in view of \re{3.13} and the fact that the slices $M_\tau$ have mean curvature
$\tau$.

From this inequality we immediate deduce
\begin{equation}
\f(t)\ge \f(t_0) e^{\frac1n (t-t_0)}\qq\A\,t\ge t_0
\end{equation}
proving the lemma.
\ep

\bl\lal{3.2}
Let $\so$ be compact and $f\in C^1(J\times \so)$, where $J$ is any open interval,
then
\begin{equation}
\f(t)=\inf_\so f(t,\cdot)
\end{equation}
is Lipschitz continuous and there holds \aev
\begin{equation}
\dot\f=\pde ft(t,x_t),
\end{equation}
where $x_t$ is a point in which the infimum is attained.

A corresponding result is also valid if $\f$ is defined by taking the supremum instead
of the infimum.
\el

\bp
$\f$ is obviously Lipschitz continuous and thus \aev differentiable by Rademacher's
theorem.

For arbitrary $t_1,t_2\in J$ we have
\begin{equation}
\f(t_1)-\f(t_2)=f(t_1,x_{t_1})-f(t_2,x_{t_2})\ge f(t_1,x_{t_1})-f(t_2,x_{t_1}).
\end{equation}

Now, let $\f$ be differentiable in $t_1$, then, by choosing $t_2>t_1$, and looking
at the difference quotients of both sides, we conclude
\begin{equation}
\dot\f(t_1)\le \pde ft(t_1,x_{t_1}).
\end{equation}

Choosing $t_2<t_1$ we obtain the opposite inequality, completing the proof of the
lemma.
\ep

We have proved that the flow hypersurfaces run straight in the singularity, if the flow
exists for all time. However, it might happen that the flow runs into the future
singularity in finite time.

To exclude this possibility we have imposed the strong volume decay condition

\bl\lal{3.3}
Let $N$ satisfy a strong volume decay condition with respect to the future, then, for
any finite $T$, $0<T\le T^*$, the flow stays in a precompact region $\Om_T$ for
$0\le t<T$.
\el

\bp
According to \rl{1.4} we may choose a time function $x^0$ such that the relation
\re{1.24} is valid for the coordinate slices $x^0=\const$.

Let $M(t)=\graph u$ be the flow hypersurfaces, and set
\begin{equation}
\f(t)=\sup_\so u(t,\cdot).
\end{equation}
Then, similarly as in the proof of \rl{3.1}, we deduce that for \aev $t$
\begin{equation}
\dot\f= \frac1{He^\psi}\le \frac1{\bar H e^\psi}\le 1,
\end{equation}
in view of \re{1.24}.

Hence we infer
\begin{equation}
\f\le \f(0)+t\qq\A\,0\le t<T^*,
\end{equation}
which proves the lemma since the  singularity corresponds to $x^0=\un$.
\ep

\section{$C^1$-estimates}

We consider a smooth solution of the evolution equation \re{0.2} in a maximal time
interval $[0,T^*)$ and shall prove a priori estimates for
\begin{equation}
\tilde v=v^{-1}=\frac1{\sqrt{1-\abs{Du}^2}}
\end{equation}
in $Q_T=[0,T]\ti \so$ for any $0<T<T^*$.

The proof is a slight modification of the proof of the corresponding result for the mean
curvature flow in \cite{cg:mz}. We note that the timelike convergence condition is
not necessary for this estimate.

Let us first state an evolution equation for $\tilde v$.

\bl[Evolution of $\tilde v$]\lal{4.1}
The quantity $\tilde v$ satisfies the evolution equation
\begin{equation}\lae{4.2}
\begin{aligned}
\dot{\tilde v}-H^{-2}\D\tilde v=&-H^{-2}\norm A^2\tilde v
-2 H^{-1}\h_{\al\bet}\n^\al\n^\bet\\
&-2H^{-2}h^{ij} x_i^\al x_j^\bet \h_{\al\bet}-H^{-2}g^{ij}\h_{\al\bet\ga}x_i^\bet
x_j^\ga\n^\al\\
&-H^{-2}\bar R_{\al\bet}\n^\al x_k^\bet\h_\ga x_l^\ga g^{kl},
\end{aligned}
\end{equation}
where $\h$ is the covariant vector field $(\h_\al)=e^{\psi}(-1,0,\dotsc,0)$.
\el

\bp
We have $\tilde v=\spd \h\n$. Let $(\x^i)$ be local coordinates for $M(t)$.
Differentiating $\tilde v$ yields the result, \cf \cite[Lemma 3.2]{cg:mz} for details.
\ep

\bl\lal{4.2}
Consider the flow in a precompact region $\Om$, then there exists a constant
$c=c(\Om)$ such that for any positive function
$0<\e=\e(x)$ on
$\so$ and any hypersurface $M(t)\su \Om$ of the flow we have
%\begin{equation}
\begin{align}
\nnorm \n&\le c\tilde v,\\\lae{4.3}
g^{ij}&\le c\tilde v^2\s^{ij},\\
\intertext{and}\lae{4.4}
\abs{h^{ij}\h_{\al\bet}x_i^\al x_j^\bet}&\le \frac{\e}{2}\norm A^2\tilde
v+\frac{c}{2\e}\tilde v^3
\end{align}
%\end{equation}
where $(\h_\al)$ is the vector field in \rl{4.1}.
\el

Confer \cite[Lemma 3.3]{cg:mz} for a proof.

Combining the preceding lemmata we infer

\bl\lal{4.3}
Consider the flow in a precompact region $\Om$, then there exists a constant
$c=c(\Om)$ such that for any positive function $\e=\e(x)$ on
$\so$ the term $\tilde v$ satisfies a parabolic inequality of the form
\begin{equation}\lae{4.6}
\dot{\tilde v}-H^{-2}\D\tilde v\le -(1-\e)H^{-2}\norm A^2\tilde
v+cH^{-2}[1+\e^{-1}]\tilde v^3.
\end{equation}
\el

\bp
The terms  on the right-hand side of \re{4.2} having a factor $H^{-2}$ can
obviously be estimated as claimed.

The remaining term can be estimated by
\begin{equation}
\begin{aligned}
2 H^{-1}\abs{\h_{\al\bet}\n^\al\n^\bet}&\le  2 c H^{-1}\tilde v^2\\
&\le \tfrac\e 2\tfrac1n \tilde v+2n c^2 \e^{-1}H^{-2}\tilde v^3.
\end{aligned}
\end{equation}

The claim then follows from the relation
\begin{equation}
\tfrac1n H^2\le \abs{A}^2,
\end{equation}
i.e.,
\begin{equation}
-H^{-2}\abs A^2\tilde v\le -\tfrac1n \tilde v.
\end{equation}
\ep

We further need the following two lemmata

\bl\lal{4.4}
Let $M(t)=\graph u(t)$ be the flow hypersurfaces, then we have
\begin{equation}
\begin{aligned}
\dot u-H^{-2}\D u&=2e^{-\psi}\tilde v H^{-1}-H^{-2}e^{-\psi}g^{ij}\bar
h_{ij}\\[\cma]
&+H^{-2}\cha 000\norm{Du}^2+2H^{-2}\cha 0i0 u^i,
\end{aligned}
\end{equation}
where the time derivative is a total derivative.
\el

\bp
We use the relation
\begin{equation}
\dot u=e^{-\psi}\tilde vH^{-1}
\end{equation}
together with \re{1.16}.
\ep

\bl\lal{4.5}
Let $\Om\su N$ be precompact and $M\su\Om$ be a spacelike graph over $\so$,
$M=\graph u$, then
\begin{equation}
\abs{\tilde v_i u^i}\le c\tilde v^3+\norm A e^\psi\norm {Du}^2,
\end{equation}
where $c=c(\Om)$.
\el

\bp
Confer the proof of \cite[Lemma 3.6]{cg:mz}.
\ep

We are now ready to prove the a priori estimate for $\tilde v$.

\bl\lal{4.6}
Let $\Om\su N$ be precompact. Then, as long as the flow stays in $\Om$, the term
$\tilde v$ is a priori bounded
\begin{equation}\lae{4.13}
\tilde v\le c=c(\Om,\sup_{M_0}\tilde v).
\end{equation}
In particular, we do not have to assume that the timelike convergence is valid, and we
note that $c$ does not depend explicitly on $T$.
\el

\bp
Let $\m,\lam$ be positive constants, where $\m$ is supposed to be small and $\lam$
large, and define
\begin{equation}\lae{4.14}
\f=e^{\m e^{-\lam u}},
\end{equation}
where we assume without loss of generality that $u\le -1$, otherwise replace in
\re{4.14} $u$ by $(u-c)$, $c$ large enough.

We shall show that
\begin{equation}
w=\tilde v \f
\end{equation}
is a priori bounded as indicated in \re{4.13} if $\m,\lam$ are chosen appropriately.

In view of \rl{4.2} and \rl{4.4} we have
\begin{equation}
\dot\f-H^{-2}\D\f\le c\m\lam e^{-\lam u}H^{-2}\tilde v^2
\f-\m\lam^2 e^{-\lam u} [1+\m e^{-\lam u}]H^{-2}\norm{Du}^2\f,
\end{equation}
since $0<H$, from which we further deduce, taking \rl{4.3} and \rl{4.5} into account,
\begin{equation}
\begin{aligned}
\dot w-H^{-2}\D w&\le -(1-\e)H^{-2} \norm A^2\tilde v\f 
+cH^{-2}[1+\e^{-1}]\tilde v^3\f\\[\cma]
&\q\,-\m\lam^2 e^{-\lam u} [1+\m e^{-\lam u}]
H^{-2}\tilde v
\norm{Du}^2\f\\[\cma]
&\q\,+c\m\lam e^{-\lam u}H^{-2}\tilde v^3\f+2\m\lam e^{-\lam u}H^{-2} \norm A
e^\psi
\norm{Du}^2\f.
\end{aligned}
\end{equation}

We estimate the last term on the right-hand side by

\begin{equation}
\begin{aligned}
2\m\lam e^{-\lam u}H^{-2}\norm A e^\psi\norm{Du}^2\f&\le (1-\e)H^{-2}\norm
A^2\tilde v\f\\ 
&+\frac{1}{1-\e}\m^2\lam^2e^{-2\lam u}H^{-2}\tilde
v^{-1}e^{2\psi}\norm{Du}^4\f,
\end{aligned}
\end{equation}
and conclude

\begin{equation}
\begin{aligned}
\dot w-H^{-2}\D w&\le  
c[1+\e^{-1}]H^{-2}\tilde v^3\f \\[\cma]
&\q\,+[\frac{1}{1-\e}-1]\m^2\lam^2 e^{-2\lam
u}H^{-2}\norm{Du}^2\tilde v\f\\[\cma]
&\q\,-\m\lam^2 e^{-\lam u}H^{-2}\norm{Du}^2\tilde v\f,
\end{aligned}
\end{equation}
where we have used that
\begin{equation}
e^{2\psi}\norm{Du}^2\le \tilde v^2.
\end{equation}

Setting $\e=e^{\lam u}$, we then obtain

\begin{equation}\lae{4.21}
\begin{aligned}
H^2(\dot w-H^{-2}\D w)&\le c e^{-\lam u} \tilde
v^3\f
+c\m\lam e^{-\lam u}\tilde v^3\f\\
&\q\,+[\frac{\m}{1-\e}-1]\m\lam^2 e^{-\lam u}\norm{Du}^2\tilde v\f.
\end{aligned}
\end{equation}

Now, we choose $\m=\frac{1}{2}$ and $\lam_0$ so large that
\begin{equation}
\frac{\m}{1-e^{\lam u}}\le \frac{3}{4}\qq\A\,\lam\ge \lam_0,
\end{equation}
and infer that the last term on the right-hand side of \re{4.21} is less than
\begin{equation}
-\frac{1}{8}\lam^2e^{-\lam u}\norm{Du}^2\tilde v\f
\end{equation}
which in turn can be estimated from above by
\begin{equation}
-c\lam^2e^{-\lam u}\tilde v^3\f
\end{equation}
at points where $\tilde v\ge 2$.

Thus, we conclude that for
\begin{equation}
\lam\ge \max (\lam_0, 4)
\end{equation}
the parabolic maximum principle, applied to $w$, yields
\begin{equation}
w\le \const (\abs{w(0)}_{_\so},\lam_0,\Om).\qedhere
\end{equation}
\ep

\section{$C^2$-estimates}\las{5}

We want to prove that, as long as the flow stays in a precompact set $\Om\su N$,
the principal curvatures of the flow hypersurfaces are a priori bounded by a constant
depending only on $\Om$ and the initial hypersurface $M_0$. Again we do not need
the timelike convergence condition for this estimate.

\cvm
Let us first prove an a priori estimate for $H$.

\bl\lal{5.1}
Let $\Om\su N$ be precompact and assume that the flow \re{0.2} stays in $\Om$
for $0\le t\le T<T^*$, then the mean curvature of the flow hypersurfaces is bounded
by
\begin{equation}
0<H\le c(\Om,\sup_{M_0}H).
\end{equation}
\el

\bp
From \rl{2.2} we immediately deduce that $\f=\log H$ satisfies the evolution
equation
\begin{equation}
\dot\f -H^{-2}\D \f=H^{-2}(\abs{A}^2+\bar
R_{\al\bet}\nu^\al\nu^\bet)-H^{-2}\norm{D\f}^2.
\end{equation}

Let $\lam$ be large and set
\begin{equation}
w=\f +\lam\tilde v.
\end{equation}
Then we conclude from \re{4.6} that $w$ satisfies the parabolic inequality
\begin{equation}
\dot w -H^{-2}\D w\le -\tfrac\lam 2 H^{-2}\abs A^2+c\lam H^{-2},
\end{equation}
if $\lam$ is large enough, $\lam\ge \lam(\Om)$. Hence the parabolic maximum
principle yields the result in view of the relation
\begin{equation}
\tfrac1n H^2\le \abs A^2.\qedhere
\end{equation}
\ep

\bl
Under the assumptions of \rl{5.1} the principal curvatures $\ka_i$, $1\le i\le n$, of
the flow hypersurfaces are a priori bounded in $\Om$
\begin{equation}
\abs{\ka_i}\le c(\Om, \sup_{M_0}\abs A).
\end{equation}
\el

\bp
Since $0\le H$, it suffices to estimate
\begin{equation}
\sup_i\ka_i\le c(\Om,\sup_{M_0}\abs A).
\end{equation}

Let $\f$  be defined  by
\begin{equation}\lae{5.8}
\f=\sup\set{{h_{ij}\h^i\h^j}}{{\norm\h=1}}.
\end{equation}
We claim that
$\f$ is a priori bounded in $\Om$.

Let $0<T<T^*$, and $x_0=x_0(t_0)$, with $ 0<t_0\le T$, be a point in $M(t_0)$ 
such that
\begin{equation}
\sup_{M_0}\f<\sup\set {\sup_{M(t)} \f}{0<t\le T}=\f(x_0).
\end{equation}

We then introduce a Riemannian normal coordinate system $(\x^i)$ at $x_0\in
M(t_0)$ such that at $x_0=x(t_0,\x_0)$ we have
\begin{equation}
g_{ij}=\de_{ij}\q \tup{and}\q \f=h_n^n.
\end{equation}

Let $\tilde \h=(\tilde \h^i)$ be the contravariant vector field defined by
\begin{equation}
\tilde \h=(0,\dotsc,0,1),
\end{equation}
and set
\begin{equation}
\tilde \f=\frac{h_{ij}\tilde \h^i\tilde \h^j}{g_{ij}\tilde \h^i\tilde \h^j}\raise 2pt
\hbox{.}
\end{equation}

$\tilde \f$ is well defined in neighbourhood of $(t_0,\x_0)$, and $\tilde \f$
assumes its maximum at $(t_0,\x_0)$. Moreover, at $(t_0,\x_0)$ we have
\begin{equation}
\dot{\tilde \f}=\dot h_n^n,
\end{equation}
and the spatial derivatives do also coincide; in short, at $(t_0,\x_0)$ $\tilde \f$
satisfies the same differential equation \re{2.8} as $h_n^n$. For the sake of
greater clarity, let us therefore treat $h_n^n$ like a scalar and pretend that
$\f=h_n^n$. 

At $(t_0,\x_0)$ we have $\dot\f\ge 0$, and, in view of the maximum principle, we
deduce from \rl{2.3}
\begin{equation}
0\le H^{-2}( -\norm A^2h_n^n + c\abs{h^n_n}^2+c).
\end{equation}

Thus $\f$ is apriori bounded in $\Om$ by a conctant $c$ depending only on $\Om$
and the initial hypersurface $M_0$.
\ep

\section{Longtime existence}

Let us look at the scalar version of the flow as in \re{3.18}
\begin{equation}\lae{6.1}
\pde ut=e^{-\psi} v H^{-1}
\end{equation}
defined in the cylinder
\begin{equation}
Q_{T^*}=[0,T^*)\times \so
\end{equation}
with initial value $u(0)\in C^\un(\so)$.

Suppose that $T^*<\un$, then, from \rl{3.3}, we conclude that the flow stays in a
compact region of $N$. Furthermore, in view of \rl{4.6} and the
$C^2$-estimates of \rs{5}, we obtain uniform $C^2$-estimates for $u$.

Thus, the differential operator on the right-hand side of \re{6.1} is uniformly elliptic in
$u$ independent of $t$, since there are constants $c_1, c_2$ such that
\begin{equation}
0<c_1\le H\le c_2\qq\A\, 0\le t<T^*,
\end{equation}
in view of \rl{2.4}.

\cvm
Hence, we can apply the known regularity results, \cf e.g., \cite[Chap. 5.5]{nk}, to
conclude that uniform $C^{2,\al}$-estimates are valid, leading further to uniform
$C^{m,\al}$-estimates for any $m\in \N$, due to the regularity result for linear
operators. But this will contradict the maximality of $T^*$.

\cvm
Therefore, $T^*=\un$, i.e., the flow exists for all time, and for any finite $T$ we have
a priori estimates in $C^m([0,T]\ti \so)$ for any $m\in \N$.

\section{A new time function}

We know that the flow exists for all time and hence we conclude from \rl{3.1} and
\rl{3.3} that the flow hypersurfaces provide a foliation of the future of $M_0$, i.e.,
the flow parameter $t$ could be used as a new time function in $D^+(M_0)$, if $Dt$
is timelike.

\bl
The flow parameter $t$ can be used as future directed time function in $D^+(M_0)$.
\el

\bp
Let $(x^\al)$ be a future directed coordinate system such that the relation \re{0.4}
is valid. Then look at the scalar version of the flow, equation \re{6.1}. If we can show
that $(\tilde x^\al)$ with
\begin{equation}
\tilde x^0=t,\qq \tilde x^i=x^i
\end{equation}
represents a regular coordinate transformation with positive Jacobi determinant,
then the lemma is proved.

\cvm
Now, the inverse coordinate transformation $x=x(\tilde x)$, which exists, since we
already know that the flow hypersurfaces provide a foliation, has the form
\begin{equation}
x^0=u(t,x)\equiv u(\tilde x),\qq x^i=\tilde x^i,
\end{equation}
where we apologize for using the same symbol $x$ to represent an $(n+1)$-tupel as
well as  the space coordinates $(x^i)$.

We immediately deduce
\begin{equation}
\Big |\pde x{\tilde x}\Big |=\pde ut >0,
\end{equation}
hence the result in view of the inverse function theorem.
\ep

The strong volume decay condition is not only sufficient to prove the long time
existence of the inverse mean curvature flow, but also necessary.

\bpp
Let $N$ be a cosmological spacetime, $M_0\su N$ a compact, spacelike hypersurface
with positive mean curvature, and suppose that the inverse mean curvature flow with
initial hypersurface $M_0$ exists for all time and provides a foliation of $D^+(M_0)$,
then $N$ satisfies a future strong volume decay condition as well as a future mean
curvature barrier condition.
\epp

\bp
Choose $x^0=t$ as new time function and let the metric of $N$ be expressed as
\begin{equation}
d\bar s^2=e^{2\psi}(-(dx^0)^2+\s_{ij}(x^0,x)dx^idx^j).
\end{equation}
$M_0$ now replaces the Cauchy hypersurface $\so$ and the flow hypersurfaces
$M(t)$ are given as graphs of functions $u$ with
\begin{equation}
u(t,x)=t.
\end{equation}

Thus we conclude from \re{6.1} that
\begin{equation}
1=\pde ut =e^{-\psi}H^{-1},
\end{equation}
or equivalently,
\begin{equation}
He^\psi=1\qq\A \,x\in M(t),
\end{equation}
i.e., the strong volume decay condition is satisfied.

\cvm
The mean curvature of the leaves $M(t)$ tends to $\un$ in view of \rl{2.4}, hence
$N$ satisfies a future mean curvature barrier condition.
\ep

From now on, let us assume that $x^0=t$ is the time function. Set
\begin{equation}
\tau=1-e^{-\frac1n t},
\end{equation}
then the future spacetime singularity corresponds to $\tau=1$, and there holds

\bt
The quantity $1-\tau$ can be looked at as the radius of the slices $\tau=\const$ as
well as a measure of the remaining life span of the spacetime, since we have
\begin{equation}
\abs{M(\tau)}=\abs{M_0} (1-\tau)^n,
\end{equation}
and the length $L(\ga)$ of any future directed curve starting from $M(\tau)$ is
estimated from above by
\begin{equation}\lae{7.10}
L(\ga)\le c (1-\tau),
\end{equation}
where
\begin{equation}\lae{7.11}
c=\frac n{\inf_{M_0}H}.
\end{equation}
\et

\bp
Let $g=\det (g_{ij})$, where $(g_{ij})$ is the induced metric of $M(t)\equiv
M(\tau)$, then
\begin{equation}
\frac d{dt}\sqrt g = -\sqrt g
\end{equation}
in view of \re{2.1}, and hence
\begin{equation}
\abs{M(t)}=\abs{M_0} e^{-t}=\abs{M_0}(1-\tau)^n.
\end{equation}

\cvm
To prove \re{7.10}, we first note that in view of \rl{2.4}
\begin{equation}
H\ge \inf_{M_0} H e^{\frac1n t}=\tfrac nc (1-\tau)^{-1},
\end{equation}
where $c$ is the constant in \re{7.11}. One of Hawking's singularity theorems then
asserts that
\begin{equation}
L(\ga)\le c (1-\tau),
\end{equation}
\cf \cite[Prop. 37 on p. 288]{bn}.
\ep
\nocite{hm}

\bibliographystyle{amsplain}
%\bibliography{mrabbrev,publications}

\providecommand{\bysame}{\leavevmode\hbox to3em{\hrulefill}\thinspace}
\providecommand{\href}[2]{#2}

%\listoffigures

%\cleardoublepage

%\thispagestyle{empty}
%\closegraphsfile
\end{document}